\documentclass[a4paper,fleqn]{cas-sc}

\usepackage[authoryear,longnamesfirst]{natbib}
\usepackage{graphicx}
\usepackage{amsmath, amssymb}
\usepackage{amsthm} 
\usepackage{appendix}
\usepackage{multirow}%
\usepackage{manyfoot}%
\usepackage{booktabs}%
\usepackage{algorithm} 
\usepackage{algorithmicx}%
\usepackage{algpseudocode}%
\usepackage{caption}
\usepackage{longtable}
\def\tsc#1{\csdef{#1}{\textsc{\lowercase{#1}}\xspace}}
\tsc{WGM}
\tsc{QE}

\newtheorem{theorem}{Theorem}
\newtheorem{lemma}[theorem]{Lemma}
\newproof{pf}{Proof}
\newtheorem{corollary}{Corollary}

\begin{document}
\let\WriteBookmarks\relax
\def\floatpagepagefraction{1}
\def\textpagefraction{.001}

\shorttitle{New Lagrangian Dual Algorithms for Solving the Continuous Nonlinear Resource Allocation Problem}    

\shortauthors{K. Hu, C. Kou and J. Yuan}  

\title [mode = title]{New Lagrangian dual algorithms for solving the continuous nonlinear resource allocation problem}  



%

\author[1,2]{Kaixiang Hu}
%
%

\address[1]{{Key Laboratory of Mathematics and Information Networks (Beijing University of Posts and Telecommunications), Ministry of Education},
	{Beijing},
	{100876}, 
	{China}}

\author[1,2]{Caixia Kou}
\cormark[1]



\address[2]{{School of science, Beijing University of Posts and Telecommunications},
	{Beijing},
	{100876}, 
	{China}}

\cortext[1]{Corresponding author}

\author[1,2]{Jianhua Yuan}
%
%


\begin{abstract}
The continuous nonlinear resource allocation problem (CONRAP) has broad applications in economics, engineering, production and inventory management, and often serves as a subproblem in complex programming. Without relying on monotonicity assumptions for the objective and constraint functions, we propose two Lagrangian dual algorithms for solving two types of CONRAP. Both algorithms determine an update strategy for the Lagrange multiplier, utilizing the values of the objective and constraint functions at the current and previous iterations. This strategy accelerates the process of finding dual optimal solutions. Subsequently, leveraging the problem's convexity, the primal optimal solution is either directly identified or derived by solving a one-dimensional linear equation. We also prove that both algorithms converge to optimal solutions within a finite number of iterations. Numerical experiments on six types of practical test problems illustrate the superior computational efficiency of the proposed algorithms. For test problems with a general inequality constraint, the first algorithm achieves a CPU time reduction exceeding an order of magnitude compared to solvers such as Gurobi and CVX. For test problems with a linear equality constraint, the second algorithm consistently outperforms four existing algorithms, delivering an improvement of over two orders of magnitude in computational efficiency.
\end{abstract}



\begin{keywords}
 Nonlinear programming\sep Large scale optimization \sep  Resource allocation problem \sep Lagrangian dual algorithm
\end{keywords}

\maketitle

\section{Introduction}\label{intro}
{In this paper, we consider the \textbf{co}ntinuous \textbf{n}onlinear \textbf{r}esource \textbf{a}llocation \textbf{p}roblem (CONRAP) with a general inequality constraint}:
\begin{align}
	\underset{\boldsymbol{x}}{\operatorname{min}} \quad &\phi(\boldsymbol{x}):=\sum_{i=1}^{n} \phi_{i}\left(x_{i}\right) \notag \\
	s.t.\quad &g(\boldsymbol{x}):=\sum_{i=1}^{n} g_{i}\left(x_{i}\right) \leq b , \label{pro1} \\ 
	&x_{i} \in X_{i}:=\left[l_{i}, u_{i}\right], \quad i=1, \ldots, n ,\notag
\end{align}
where $ \phi_{i}: \mathbb{R} \rightarrow \mathbb{R} \text { and } g_{i}: \mathbb{R} \rightarrow \mathbb{R} $ are convex and continuously differentiable, $ b \in \mathbb{R} $. In the following content, we also use $ \boldsymbol{x} \in \boldsymbol{X} $ to represent the box constraints.

{We also consider the CONRAP with a linear equality constraint}:
\begin{align}
	\underset{\boldsymbol{x}}{\operatorname{min}} \quad &\phi(\boldsymbol{x}) =\sum_{i=1}^{n} \phi_{i}\left(x_{i}\right) \notag \\
	s.t.\quad &g(\boldsymbol{x}):=\sum_{i=1}^{n} a_{i} x_{i} = b, \label{pro2}  \\ 
	&x_{i} \in \left[l_{i}, u_{i}\right], \quad i=1, \ldots, n,  \notag
\end{align}
where $ a_{i} \neq 0 $, and the sign is all the same for $ i=1, \ldots, n $. Both problems are convex, indicating that this paper focuses solely on the convex types of CONRAP. 
\subsection{Applications}
{The CONRAP arises in various fields, including the resource allocation problem in commodity warehousing, portfolio investment, statistics, as well as the balance problem in engineering and economics} (\citet{patriksson2008survey}, \citet{maloney1993constrained}, \citet{bretthauer1999nonlinear}, \citet{dussault1986convex}). Additionally, this problem serves as a subproblem in algorithms proposed for complex programming problems (\citet{patriksson2015algorithms}, \citet{cipolla2022training}).
Five types of practical problems are given as follows.

{\textbf{Commodity warehousing problem}. {This problem aims to achieve a balance between the holding cost and order quantity while satisfying the storage capacity constraint} (\citet{nielsen1993massively}, \citet{katoh1998resource}, \citet{gao2017economic}, \citet{dai2018location}, \citet{kamaludin2021agricultural}).  The initial problem model is as follows:
	\begin{align*}
		\underset{\boldsymbol{x}}{\operatorname{min}} \quad &\sum_{i=1}^{n} c_{i}x_{i} + k_{i}/x_{i}  \\
		s.t.\quad &\sum_{i=1}^{n} a_{i}x_{i} \leq b , \\ 
		&x_{i} =\left[l_{i}, u_{i}\right], \quad i=1, \ldots, n ,
	\end{align*}
	where $ x_{i} $ denote the order quantity of item $ i $, $ c_{i} $ is the holding cost, $ k_{i} $ is the ordering (or replenishment) cost, $ a_{i} $ is the storage requirement per item and $ b $ is the storage capacity.
	
	\textbf{Euclidean projection problem}. As a kind of application of problem (\ref{pro2}), {Euclidean projection problems} arise in many optimization problems, especially as subproblems of complex optimization problems (\citet{ventura1988note}, \citet{maloney1993constrained}, \citet{dai2006new}, \citet{dattorro2010convex}, \citet{gabidullina2018problem}, \citet{usmanova2021fast}). The problem model is as follows:
	\begin{align*}
		\underset{\boldsymbol{x}}{\operatorname{min}} \quad &\sum_{i=1}^{n} \frac{1}{2} (x_{i}-y_{i})^2 \\
		s.t.\quad &\sum_{i=1}^{n} a_{i}x_{i} = b , \\ 
		&x_{i} =\left[l_{i}, u_{i}\right], \quad i=1, \ldots, n ,
	\end{align*}
	where $ \boldsymbol{x} $  denotes the projection of any given vector $ \boldsymbol{y} \in \mathbb{R}^n$.
	
	\textbf{Portfolio investment problem}. {This problem aims to find the optimal investment plan under a fixed total asset}, which can be found in \citet{nielsen1992massively}, \citet{kiwiel2007linear}, \citet{kiwiel2008variable}, \citet{kiwiel2008breakpoint} and  improved by \citet{chen2015artificial}, \citet{seyedhosseini2016novel}, \citet{li2019high}, \citet{cura2021rapidly}. The problem model is as follows:
	\begin{align*}
		\underset{\boldsymbol{x}}{\operatorname{min}}\quad&\sum_{i=1}^{n} \frac{1}{2}d_{i} x_{i}^2 - c_{i} x_{i} \\
		s.t.\quad&\sum_{i=1}^{n} x_{i} = 1,  \\
		&x_{i} \in \left[l_{i}, u_{i}\right], \quad i=1, \ldots, n, 
	\end{align*}
	where $ x_{i} $  denotes the investment proportion of product $ i $, $ d_{i} $ represents a diagonal approximation of the positive definite covariance matrix, $ c_{i} $ is the expected asset returns.
	
	\textbf{Sampling problem}. {This problem aims to find the optimal stratified sampling strategy under a fixed total sample size} 
 (\citet{bretthauer1999nonlinear}, \citet{khan2015designing}, \citet{shields2015refined}, \citet{varshney2017optimum}, \citet{varshney2019optimum}, \citet{nguyen2021stratified}, \citet{alshqaq2022nonlinear}). The problem model is as follows:
	\begin{align*}
		\underset{\boldsymbol{x}}{\operatorname{min}}\quad &\sum_{i=1}^{n} \omega_{i}^{2} \frac{\left(M_{i}-x_{i}\right) \sigma_{i}^{2}}{\left(M_{i}-1\right) x_{i}}  \\
		s.t.\quad &\sum_{i=1}^{n} x_{i}=b, \\ 
		&x_{i} \geqslant 1, \quad i=1, \ldots, n ,
	\end{align*}
	where $ x_{i}$ is the chosen samples in each strata, $ M_{i} $ is the population in each strata, $ M $ is the entire population, $ \omega_{i}=M_{i} / M, \sigma_{i}^{2} $  is an appropriate estimate of the variance in each strata, $ b $ is the total sample size, and at least one sample is taken from each strata. 
	
	\textbf{Target search problem}. {This problem aims to obtain the optimal search strategy under a fixed search time}, which is firstly described in \citet{koopman1953optimum}, and then studied by \citet{stone1981search}, \citet{koopman1999search}, \citet{hohzaki2006search}, \citet{chen2010tracking}, \citet{hohzaki2015search}, \citet{vaillaud2023target}. The problem model is as follows:
	\begin{align*}
		\underset{\boldsymbol{x}}{\operatorname{min}}\quad &\sum_{i=1}^{n}m_{i}\left(e^{-c_{i} x_{i}}-1\right) \\
		s.t.\quad &\sum_{i=1}^{n} x_{i}=b, \\ 
		&x_{i} \in \left[l_{i}, u_{i}\right], \quad i=1, \ldots, n,
	\end{align*}
	where $ x_{i}$ is the search time for the box $ i $, $ m_{i} $ is the probability of an target being inside box $ i $, and $ -c_{i} $ is proportional to the difficulty of searching inside the box $ i $. The objective is to maximize the overall probability of finding the target, $\sum_{i=1}^{n} -m_{i}\left(e^{-c_{i} x_{i}}-1\right) $.
}
\subsection{Related work}
There has been extensive research on the CONRAP, resulting in the development of numerous efficient algorithms. {The resource allocation optimization problem in search theory is formulated as a special problem (\ref{pro1}) without box constraints (\citet{charnes1958theory}).  The first Lagrangian dual algorithm for problem (\ref{pro1}) is developed, where the Lagrange multiplier  is updated by the bisection method.} Additionally, a pioneering study by \citet{luss1975allocation} formulates the resource allocation problem among multiple competitive activities as a special problem (\ref{pro2}) without box constraints. {The first pegging algorithm is proposed by making} the most of this problem's separability.

The pegging algorithm is a specialized approach developed for the CONRAP, exploiting the problem's separability. At each iteration, the algorithm fixes certain variables to satisfy the box constraints, with the pegging rules derived from the KKT conditions. Through the iteration, this algorithm converges to the optimal solution. This approach is further developed by \citet{bitran1977design}, who propose the classical pegging algorithm. Nevertheless, this algorithm is unable to solve problem (\ref{pro1}) featuring the nonlinear constraint. To address this issue, Bretthauer and Shetty (\citet{bretthauer2002pegging}, \citet{bretthauer2002nonlinear}) propose the modified pegging algorithm, designed specifically for CONRAP with monotonic objective and constraint functions. The core idea involves relaxing the box constraints and solving a series of subproblems. Variables violating the box constraints are identified and updated to their respective bounds. This method is later refined into the 3-sets pegging algorithm (\citet{kiwiel2008breakpoint}) and the 5-sets pegging algorithm (\citet{de2012breakpoint}), enhancing computational efficiency for large-scale problems.

{Meanwhile, research on algorithms for CONRAP can be broadly divided into two categories based on differences in the objective function.}  
The CONRAP featuring a quadratic objective function has been widely studied in the past decades. Importantly, \citet{dai2006new} propose a novel Lagrangian dual algorithm with linear time complexity for separable \textbf{s}ingly \textbf{l}inearly constrained \textbf{q}uadratic
\textbf{p}rograms subject to lower and upper \textbf{b}ounds (SLBQP), which  is referred to as the SLBQP1 algorithm in this paper. In the SLBQP1 algorithm, the Lagrange multiplier is updated using an improved secant-like method to ensure the solution satisfies the linear constraint $ g(\boldsymbol{x}) = b $. Notably, this algorithm is specifically designed for CONRAP with a quadratic objective function (e.g., the portfolio investment problem), as it relies on the explicit solution to the Lagrangian dual problem. 
Additionally, a projection method is employed to enforce the box constraints $ \boldsymbol{x} \in \boldsymbol{X} $.

Furthermore, the optimal solution to the quadratic CONRAP can be obtained by selecting the median of the active boundary points and iteratively shrinking the boundary point set (\citet{kiwiel2008breakpoint}). The necessary and sufficient condition for the optimal solution is established based on the KKT conditions, and the optimal solution is identified using the bisection method (\citet{pardalos1990algorithm}).
Moreover, a lightweight open-source software library is developed for this type of problems (\citet{frangioni2013library}).

For the non-quadratic CONRAP, the theoretical properties of problem (\ref{pro2}) are evaluated, and a necessary and sufficient condition for identifying the optimal solution is established based on the KKT conditions (\citet{stefanov2015solution}). Additionally, a penalized algorithm based on Bregman distance is proposed with the explicit projection, and this algorithm transforms the problem into a minimization problem under a unit simplex constraint (\citet{hoto2020penalty}). 
Moreover, an efficient algorithm based on the augmented Lagrangian method is proposed by \citet{torrealba2022augmented} and is referred to as the Aug\_lag algorithm in this paper. It reformulates the original problem into a sequence of simpler subproblems, solved by using the Newton method. To enforce constraints, it augments the Lagrangian dual function with a penalty term.
This algorithm alternates between solving the augmented Lagrangian function to update the solution and projecting the solution onto the box constraints.

However, taken together, these studies indicate that most existing algorithms impose restrictions on the linearity of constraints or the monotonicity of both objective and constraint functions, limiting their applicability. Moreover, no prior research has integrated the efficient Lagrangian dual algorithm with the separability of the problem to solve CONRAP. Motivated by these factors, we attempt to propose  Lagrangian dual algorithms that leverage the separability of the problem to accelerate the search for the optimal solution while ensuring convergence guarantees.
\subsection{Our contributions}
The main contributions of this paper include:
\begin{itemize}
	{\item Targeting the two types of the CONRAP, namely problem (\ref{pro1}) and problem (\ref{pro2}), we propose two novel Lagrangian dual algorithms that leverage the convexity and separability of these problems. Both algorithms do not rely on monotonicity assumptions, which significantly broadens their applicability. }Based on the convexity of the objective and constraint functions, we design a Lagrange multiplier update strategy by the values of these functions. Utilizing the problem's separability, the Lagrangian dual problem is decomposed into $ n $ one-dimensional tractable subproblems; these subproblems are rapidly solved by using the function's mononicity and gradient.
	\item Through rigorous theoretical analysis, we prove that both of the proposed algorithms converge to optimal solutions within a finite number of iterations. By truncating the Lagrange multiplier using a bisection strategy, the interval containing the optimal multiplier is guaranteed to shrink by a fixed proportion. This guarantees that both algorithms terminate after a finite number of iterations.
	Subsequently, by leveraging the convexity of both objective and constraint functions, we prove that the obtained solution satisfies the KKT system of the convex problem, that is, the optimal solution is identified. 
	\item {Extensive numerical experiments demonstrate that both of the proposed algorithms consistently identify the optimal solution across all test problems as well as exhibit a significant advantage in computational efficiency. Notably, across all six types of test problems, both algorithms exhibit superior efficiency, outperforming four existing methods, including Gurobi solver and CVX, by at least an order of magnitude.} 
\end{itemize}
\subsection{Organization of the paper}
The remainder of this paper is organized as follows. Sections \ref{sec2} and \ref{sec3} propose the corresponding Lagrangian dual algorithms for the problem (\ref{pro1}) and (\ref{pro2}), respectively. Section \ref{sec4} presents the convergence analysis of the proposed algorithms. Section \ref{sec5} presents numerical experiments to evaluate the performance of the algorithms and discusses the results. Finally, Section \ref{sec6} concludes the paper.

\section{A new Lagrangian dual algorithm for problem (\ref{pro1})}  \label{sec2}
{In this section, we propose an  exact Lagrangian dual algorithm for problem (\ref{pro1}).
The Lagrangian dual problem corresponding to problem (\ref{pro1}) is}
\begin{equation}\label{lag1}
	\underset{\lambda \geq 0}{\operatorname{max}} \underset{\boldsymbol{x} \in \boldsymbol{X}}{\operatorname{min}} \quad L(\boldsymbol{x},\lambda)= \phi(\boldsymbol{x})+\lambda( g(\boldsymbol{x})-b),
\end{equation}
where $ L(\boldsymbol{x},\lambda) $ is the Lagrangian dual function of problem (\ref{pro1}).  
{For a fixed Lagrange multiplier $\lambda$, we solve the following problem to obtain the optimal solution to $  \underset{\boldsymbol{x} \in \boldsymbol{X}}{\operatorname{min}} \  L(\boldsymbol{x},\lambda) $,
\begin{equation}\label{sub1}
	\underset{\boldsymbol{x} \in \boldsymbol{X}}{\operatorname{min}} \quad \phi(\boldsymbol{x})+\lambda g(\boldsymbol{x}),
\end{equation} 
until the optimal multiplier $  \lambda^* $ is  identified. Then the optimal solution to $  \underset{\boldsymbol{x} \in \boldsymbol{X}}{\operatorname{min}} \  L(\boldsymbol{x},\lambda^*) $ yeilds the optimal solution to problem (\ref{pro1}), which can be guaranteed by the convexity of this problem. Therefore,  how to rapidly identify the optimal multiplier $ \lambda^* $ is the key to our algorithm.}

\subsection{Lagrange multiplier update strategy} \label{mum}
The Lagrange multiplier $ \lambda $ is updated based on the values of both the objective and constraint functions at current and previous iterations until the optimal multiplier $ \lambda^* $ is identified.
Note that $ \lambda $ takes values in the range $ [0, +\infty) $. When $ \lambda = 0 $ or $ \lambda = +\infty $, problem (\ref{sub1}) reduces to $ \underset{\boldsymbol{x} \in \boldsymbol{X}}{\operatorname{min}} \ \phi(\boldsymbol{x}) $ or $ \underset{\boldsymbol{x} \in \boldsymbol{X}}{\operatorname{min}} \ g(\boldsymbol{x}) $, respectively. We denote the corresponding solutions as  $ \boldsymbol{x_{\phi}} = \underset{\boldsymbol{x} \in \boldsymbol{X}}{\operatorname{argmin}\  }\phi(\boldsymbol{x}) $ and $ \boldsymbol{x_{g}} = \underset{\boldsymbol{x} \in \boldsymbol{X}}{\operatorname{argmin}\  }g(\boldsymbol{x}) $, respectively. If $ g(\boldsymbol{x_{\phi}}) \leq b $, the optimal solution is found; if $ g(\boldsymbol{x_{g}})>b $, the problem (\ref{pro1}) has no feasible solutions. Therefore, we only consider the case of $ g(\boldsymbol{x_{g}}) \leq b  $ and $ g(\boldsymbol{x_{\phi}}) > b $ in the following content. In this case, we must increase $ \lambda $ to enhance the influence of the constraint function in $ L(\boldsymbol{x},\lambda) $. 

\begin{figure}[h]
	\centering
	\includegraphics[scale=0.65]{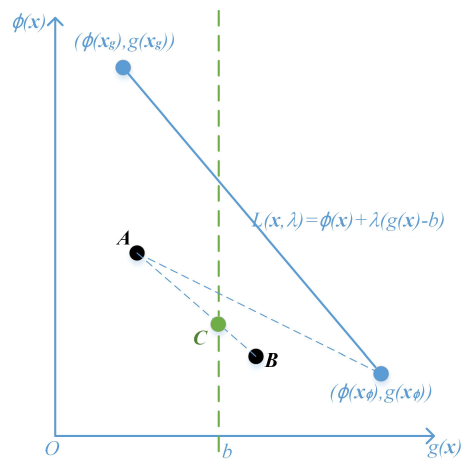}
	\caption{{In this figure}, the vertical coordinate axis consists of the function $\phi(\boldsymbol{x}) $ value while the horizontal coordinate axis consists of the function  $g(\boldsymbol{x}) $ value, the initial iterates  $ \boldsymbol{x_{g}} = \underset{\boldsymbol{x} \in \boldsymbol{X}}{\operatorname{argmin}\  }g(\boldsymbol{x}) $ and $ \boldsymbol{x_{\phi}} = \underset{\boldsymbol{x} \in \boldsymbol{X}}{\operatorname{argmin}\  }\phi(\boldsymbol{x}) $, and each iterate $\boldsymbol{x} $ fulfills the box constraints $ \boldsymbol{x} \in \boldsymbol{X} $. }
	\label{fig1}
\end{figure}
Geometrically, the value of "$ -\lambda $" represents the slope of the line $ \phi(\boldsymbol{x}) + \lambda g(\boldsymbol{x}) = \textit{constant} $ in the $ (g(\boldsymbol{x}), \phi(\boldsymbol{x})) $-plane, as illustrated in Figure \ref{fig1}. 
	For the initial iterates \( \boldsymbol{x_{g}} \) and \( \boldsymbol{x_{\phi}} \), there exists a line connecting the points \( (g(\boldsymbol{x_{g}}), \phi(\boldsymbol{x_{g}})) \) and \( (g(\boldsymbol{x_{\phi}}), \phi(\boldsymbol{x_{\phi}})) \). 
	Thus,  the initial multiplier value is set to the negative slope value of this line, given by $\lambda=-\frac{\phi\left(\boldsymbol{x_{g}}\right)-\phi\left(\boldsymbol{x_{\phi}}\right)}{g\left(\boldsymbol{x_{g}}\right)-g\left(\boldsymbol{x_{\phi}}\right)} $.  
	
	During the iteration, the new iterate \( \boldsymbol{x} \) is obtained by solving problem (\ref{sub1}). And  \( \boldsymbol{x} \) is stored as \( \boldsymbol{x}_{+} \) if \( g(\boldsymbol{x}) > b \), or as \( \boldsymbol{x}_{-} \) if \( g(\boldsymbol{x}) < b \). The new multiplier is updated as \( \lambda:= -\frac{\phi\left(\boldsymbol{x}_{-}\right) - \phi\left(\boldsymbol{x}_{+}\right)}{g\left(\boldsymbol{x}_{-}\right) - g\left(\boldsymbol{x}_{+}\right)} \), since for this multiplier, the condition $ \phi(\boldsymbol{x}_{+})+\lambda g(\boldsymbol{x}_{+}) = \phi(\boldsymbol{x}_{-})+\lambda g(\boldsymbol{x}_{-}) $ holds. 
	This process continues until either the iterate $ \boldsymbol{x} $ satisfies $ g(\boldsymbol{x}) = b $, or the multiplier value no longer changes. In these two cases, the optimal multiplier $ \lambda^* $ is identified, respectively. Details are provided below:
	\begin{itemize}
		\item[(I)] For the first case, 
		given a multiplier $\lambda$, by minimizing $L(\boldsymbol{x},\lambda)$, we obtain its optimal solution $\boldsymbol{x}$. If $g(\boldsymbol{x})=b$ is satisfied, then according to Theorem \ref{thm_op1}, the current $\lambda$ is guaranteed to be the optimal multiplier, and $\boldsymbol{x}$ is the optimal solution to problem (\ref{pro1}), as shown by point $ C $ of Figure \ref{fig1}.
		\item[(II)] For the second case, given a multiplier $\lambda$, by minimizing $L(\boldsymbol{x},\lambda)$,
		if its optimal solution  $ \boldsymbol{x} $ satisfies $ L(\boldsymbol{x},\lambda)= L(\boldsymbol{x}_{+},\lambda)=L(\boldsymbol{x}_{-},\lambda)$ for this multiplier, where $ \boldsymbol{x}_{+} , \boldsymbol{x}_{-} $ are the two preceding iterates, then according to Lemma \ref{lemma2}, this multiplier is shown to be the optimal multiplier $ \lambda^* $. 
		As shown in Figure \ref{fig1}, when the iteration process terminates, points $A$ and $B$ are obtained, but the optimal solution $C$ lies on the line connecting $A$ and $B$. Therefore, further steps are needed to identify the optimal solution $ \boldsymbol{x}^* $.
\end{itemize}
Fortunately, in the second case, the information from the last two iterations can be leveraged to facilitate the identification process. As demonstrated in Lemma \ref{lem_omega}, the optimal solution set of $ \underset{\boldsymbol{x}\in \boldsymbol{X}}{\operatorname{min}\  }L(\boldsymbol{x},\lambda^*) $, denoted as $ \Omega^*=\{\boldsymbol{x}|\boldsymbol{x}=\underset{\boldsymbol{x}\in \boldsymbol{X}}{\operatorname{argmin}\  }L(\boldsymbol{x},\lambda^*)\} $, is convex. As established in the proof of Theorem \ref{thm_op2}, the function $ g(\boldsymbol{x}) $ is linear for $ \boldsymbol{x} \in \Omega^* $. The optimal solution is identified by taking the convex combination of $\boldsymbol{x}_{+}$ and $\boldsymbol{x}_{-}$, and solving the following linear subproblem.
\begin{equation}\label{sub2}
	\{\alpha \ \mid \ g(\alpha \boldsymbol{x}_{+}+(1-\alpha)\boldsymbol{x}_{-}) - b=0 \}.
\end{equation}
The solution identified by this subproblem is proved to be the optimal solution to problem (\ref{pro1}) and the optimality proof is given as  Theorem \ref{thm_op2} in Section \ref{sec4}.
	
	Furthermore, in the above content, it is imperative to solve the  problem (\ref{sub1}) to obtain the iterates, and $ \lambda \in [0,+\infty) $. By taking advantage of the convexity and separability of Lagrangian dual function, this problem is decomposed into $ n $ one-dimensional tractable subproblems and readily solved by using the function's mononicity and gradient.

\subsection{Algorithm steps} \label{alstep}
Based on the above analysis, our first algorithm is given as Algorithm \ref{algo1}.
\begin{algorithm}[h]
	\caption{(Identify the optimal solution to problem (\ref{pro1}))} \label{algo1} 
	$ \text { Procedure }\left[ \boldsymbol{x}^*\right]=\operatorname{Algorithm\ 1}(\phi(\boldsymbol{x}),g(\boldsymbol{x}) , b, X) $ 
	\begin{algorithmic}[1]
		\State Initialization:
		set a constant $  \gamma \in(0,\dfrac{1}{2}) $ and  $\epsilon >0 $.\\ 
		$ \boldsymbol{x}_{+}:=\underset{\boldsymbol{x} \in \boldsymbol{X}}{\operatorname{argmin}}\  \phi(\boldsymbol{x}) $, and $  \underline{\theta}:=0 $.
		\If {$ g\left(\boldsymbol{x}_{+}\right) \leq b  $} 
		\State return $ \boldsymbol{x}^*= \boldsymbol{x}_{+} $.
		\Else \State $ \boldsymbol{x}_{-}:=\underset{\boldsymbol{x} \in \boldsymbol{X}}{\operatorname{argmin}}\  g(\boldsymbol{x}) $ , and $ \bar{\theta}:=\frac{\pi}{2} $.
		\EndIf
		\If {$  g\left(\boldsymbol{x}_{-}\right)>b  $} 
		\State return "There is no solution."
		\ElsIf{$  g\left(\boldsymbol{x}_{-}\right)=b  $} 
		\If {$ \|\nabla g\left(\boldsymbol{x}_{-}\right)\| \neq 0 $} 
		\State return $  \boldsymbol{x}^*= \boldsymbol{x}_{-} $. 
		\Else \State $  \boldsymbol{x}^*=\operatorname{Projected \_Gradient}(\phi(\boldsymbol{x}),g(\boldsymbol{x}) , \boldsymbol{x}_{-},  X) $ and return  $ \boldsymbol{x}^* $. 
		\EndIf
		\Else 
		\While{$ \bar{\theta}-\underline{\theta} > \epsilon $} \label{while_begin}
		\State $\lambda:=\frac{\phi\left(\boldsymbol{x}_{-}\right)-\phi\left(\boldsymbol{x}_{+}\right)}{g\left(\boldsymbol{x}_{+}\right)-g\left(\boldsymbol{x}_{-}\right)} $ and $   \theta:=\arctan \lambda  $. \label{lam_up1}
		
		\If {$ \left|\theta-\frac{\bar{\theta}+\underline{\theta}}{2}\right|>\left(\frac{1}{2}-\gamma\right)(\bar{\theta}-\underline{\theta}) $} 
		\State set  $ \theta:=\frac{\bar{\theta}+\underline{\theta}}{2} $  and $  \lambda:=\tan \theta  $.\label{lam_up2}
		\EndIf
		\State Solve problem  (\ref{sub1}) to obtain $ \tilde{\boldsymbol{x}} $.
		\If {$ g\left(\tilde{\boldsymbol{x}}\right) = b $}
		\State {return $  \boldsymbol{x}^*= \tilde{\boldsymbol{x}}$.} \label{optimum1}
		\ElsIf{$  g\left(\tilde{\boldsymbol{x}}\right)>b  $}
		\State set $  \boldsymbol{x}_{+}=\tilde{\boldsymbol{x}} $ and $  \underline{\theta}=\theta $.
		\Else \State set $  \boldsymbol{x}_{-}=\tilde{\boldsymbol{x}} $ and $   \bar{\theta}=\theta $.
		\EndIf
		\EndWhile \label{while_end}
		\EndIf
		\State Solve the linear subproblem  $  g(\alpha^* \boldsymbol{x}_{+}+(1-\alpha^*)\boldsymbol{x}_{-}) - b=0 $ to obtain $ \alpha^*  $
		\State $  \boldsymbol{x}^* = \alpha^* \boldsymbol{x}_{+}+(1-\alpha^*)\boldsymbol{x}_{-}$ and return  $ \boldsymbol{x}^* $.\label{optimum2}
	\end{algorithmic}
\end{algorithm}
Firstly, we set the interval $[\underline{\theta},\bar{\theta}] $,  $ \theta:=\arctan \lambda $ and $\theta \in [\underline{\theta},\bar{\theta}] $ always holds in this algorithm.
In order to ensure that this interval is  reduced by a certain multiple during the algorithm's iteration, we adopt the update strategy  
proposed in \citet{kou2019bisection}. 
When the condition $ \left|\theta-\frac{\bar{\theta}+\underline{\theta}}{2}\right|\leq \left(\frac{1}{2}-\gamma\right)(\bar{\theta}-\underline{\theta}) $ is satisfied ($ \gamma $ is a given parameter), we adopt the current value of multiplier $ \lambda $. Otherwise, the current $ \theta $ is too close to the interval boundary and  we set $ \theta=\frac{\bar{\theta}+\underline{\theta}}{2} $. Therefore, the interval length decreases by
at least a multiple of $ (1-\gamma) $ in each iteration, which is contributed to conduct the convergence analysis of the proposed algorithm.

{As discussed in Section \ref{mum}, we  only need to consider the case of $ g(\boldsymbol{x_{g}}) \leq b  $ and $ g(\boldsymbol{x_{\phi}}) > b $. Particularly, when $  g\left(\boldsymbol{x}_{g}\right)=b  $ holds, it is necessary to identify the optimal solution to $ \phi(\boldsymbol{x}) $ under the set $ \Omega_g $, where $ \Omega_g =\{\boldsymbol{x}\mid \boldsymbol{x}=\underset{\boldsymbol{x}\in \boldsymbol{X}}{\operatorname{argmin}\  }g(\boldsymbol{x})\} $  is the optimal solution set of $ g(\boldsymbol{x}) $ under box constraints.
	If $ \|\nabla g\left(\boldsymbol{x}_{g}\right)\| \neq 0 $, then $ \Omega_g $ has only one element $ \boldsymbol{x}_{g} $. Otherwise, $ \Omega_g $ may  have multiple elements and we adopt the  projected gradient algorithm proposed by \citet{iusem2003convergence} to identify the solution $ \boldsymbol{x} $ satisfying $ \boldsymbol{x}= \underset{\boldsymbol{x}\in \Omega_g}{\operatorname{argmin}\  }\phi(\boldsymbol{x})  $. 
	The details  of this algorithm  are presented in Appendix \ref{appa}. 
	
	When $ g(\boldsymbol{x_{g}}) < b  $ and $ g(\boldsymbol{x_{\phi}}) > b $, we update the multiplier $ \lambda $ and iterate $ \boldsymbol{x} $ by utilizing  the proposed Lagrange multiplier update strategy  and solving the problem (\ref{sub1}), respectively. Finally, the optimal solution $ \boldsymbol{x}^* $ is identified at line \ref{optimum1} or \ref{optimum2}, corresponding to the two cases described in Section \ref{mum}. 
	
	\section{A new Lagrangian dual algorithm for problem (\ref{pro2})} \label{sec3}
	{Given that problem (\ref{pro2}) shares analogous structure and properties with problem (\ref{pro1}), an exact Lagrangian dual algorithm is proposed by invoking Algorithm 1 in distinctive cases.
	Without losing generality, we assume $ a_i >0,\ \forall i=1,\ldots,n $, since the sign of $ a_i $ is all the same in the constraint $ g(\boldsymbol{x}) = \sum_{i=1}^{n} a_{i} x_{i} = b  $. This indicates that the function $ g(\boldsymbol{x}) $ is a monotone increasing linear function.}
		First of all, we present the Lagrangian dual problem of problem (\ref{pro2})
		\begin{equation}\label{lag2}
			\underset{\lambda }{\operatorname{max}} \underset{\boldsymbol{x} \in \boldsymbol{X}}{\operatorname{min}} \quad L(\boldsymbol{x},\lambda)= \phi(\boldsymbol{x})+\lambda( g(\boldsymbol{x})-b), 
		\end{equation}
		where $ L(\boldsymbol{x},\lambda) $ is the Lagrangian dual function of problem (\ref{pro2}). 

		\begin{figure}[h]
			\centering
			\includegraphics[scale=0.65]{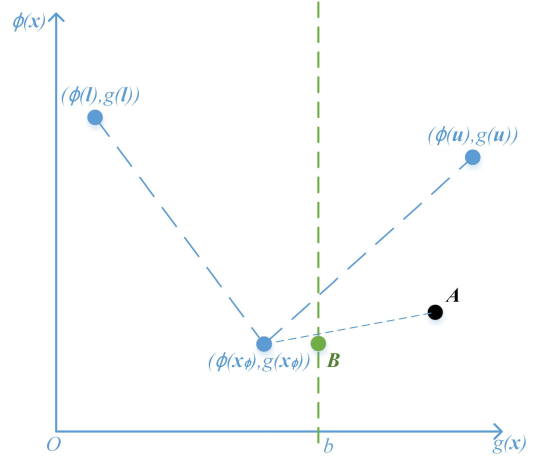}
			\caption{{In this figure}, the vertical coordinate axis consists of the function $\phi(\boldsymbol{x}) $ value while the horizontal coordinate axis consists of the function  $g(\boldsymbol{x}) $ value, and each iterate $x $ fulfills the box constraints, $ \boldsymbol{x} \in \boldsymbol{X} $. And $ \boldsymbol{x_\phi} = \underset{\boldsymbol{x} \in \boldsymbol{X}}{\operatorname{argmin}\  }\phi(\boldsymbol{x}) $. }
			\label{fig3}
		\end{figure}
		Unlike problem (\ref{pro1}), where the multiplier must be non-negative, the multiplier associated with the equality constraint in problem (\ref{pro2}) can take negative values. Consequently, the region containing the optimal solution extends to both sides of the vertical line $ g(\boldsymbol{x}) = b $, as shown in Figure \ref{fig3}.
		When $ \lambda = +\infty $ or $ \lambda = -\infty $, problem (\ref{sub1}) reduces to $ \underset{\boldsymbol{x} \in \boldsymbol{X}}{\operatorname{min}} \ g(\boldsymbol{x}) $ or $ \underset{\boldsymbol{x} \in \boldsymbol{X}}{\operatorname{max}} \ g(\boldsymbol{x}) $, respectively.  The corresponding solutions is $ \boldsymbol{l} $ and $ \boldsymbol{u} $, respectively.
		If $ g(\boldsymbol{l})>b  $ or $ g(\boldsymbol{u})<b $, the problem (\ref{pro2}) has no feasible solutions. Therefore, we only consider the case of $ g(\boldsymbol{l}) \leq b  $ and $ g(\boldsymbol{u}) \geq b $ in the following content. Particularly, when $ g(\boldsymbol{x_{\phi}}) = b $, where $ \boldsymbol{x_{\phi}} = \underset{\boldsymbol{x} \in \boldsymbol{X}}{\operatorname{argmin}\  }\phi(\boldsymbol{x}) $, the optimal solution is found.
		
		Importantly, the value of "$ -\lambda $" still corresponds to the slope of the line $ \phi(\boldsymbol{x}) + \lambda g(\boldsymbol{x}) = \textit{constant} $ in the $ (g(\boldsymbol{x}), \phi(\boldsymbol{x})) $-plane. Specifically, the multiplier is updated as: \( \lambda:= -\frac{\phi\left(\boldsymbol{x}_{-}\right) - \phi\left(\boldsymbol{x}_{+}\right)}{g\left(\boldsymbol{x}_{-}\right) - g\left(\boldsymbol{x}_{+}\right)} \). When $ g(\boldsymbol{x_{\phi}}) < b $, we initialize $ \boldsymbol{x}_{-} =\boldsymbol{x}_{\phi}$ and $\boldsymbol{x}_{+}=\boldsymbol{u} $,  then invoke Algorithm 1 to determine the optimal solution. Conversely, when $ g(\boldsymbol{x_{\phi}}) > b  $,  we initialize $ \boldsymbol{x}_{-} =\boldsymbol{l}$ and $\boldsymbol{x}_{+}=\boldsymbol{x}_{\phi} $, followed by a call to Algorithm 1.  
	

	For problem (\ref{pro2}), our algorithm is given as Algorithm \ref{algo4}.
	
	\begin{algorithm}[h] 
		\caption{(Identify the optimal solution to problem (\ref{pro2}))} \label{algo4}
		$ \text { Procedure }\left[ \boldsymbol{x}^*\right]=\operatorname{Algorithm\  2}(\phi(\boldsymbol{x}),g(\boldsymbol{x}) , b, X) $ 
		\begin{algorithmic}[1]
			\State Initialization: if $ g(\boldsymbol{l})>b  $ or $ g(\boldsymbol{u})<b $, then return "There is no solution."
			\State $ \boldsymbol{x}_{+}:=\underset{\boldsymbol{x} \in \boldsymbol{X}}{\operatorname{argmin}}\  \phi(\boldsymbol{x}) $.
			\If {$ g\left(\boldsymbol{x}_{+}\right) = b  $} 
			\State  {return $ \boldsymbol{x}^* = \boldsymbol{x}_{+} $.}
			\ElsIf{$ g\left(\boldsymbol{x}_{+}\right) < b  $}
			\State  set $ (\boldsymbol{x}_{-},\boldsymbol{x}_{+})=(\boldsymbol{x}_{+},\boldsymbol{u}) $, and $ (\bar{\theta},\underline{\theta}) = (0,-\frac{\pi}{2}) $.
			\State  $ \boldsymbol{x}^* = \operatorname{Algorithm\  1}(\phi(\boldsymbol{x}),g(\boldsymbol{x}) , b, X) $ and return  $ \boldsymbol{x}^* $.
			\Else \State set $ \boldsymbol{x}_{-}=\boldsymbol{l} $, and $ (\bar{\theta},\underline{\theta}) = (\frac{\pi}{2},0) $.
			\State $ \boldsymbol{x}^* = \operatorname{Algorithm\  1}(\phi(\boldsymbol{x}),g(\boldsymbol{x}) , b, X) $ and return  $ \boldsymbol{x}^* $.
			\EndIf
		\end{algorithmic}
	\end{algorithm}

\section{Convergence analysis}  \label{sec4}
{In this section, we prove that both of the proposed algorithms terminate in finite iterations and identify the optimal solution. For distinguished termination steps of the algorithm, we present optimality proofs of the corresponding returned solutions.}

Firstly, we proceed with the convergence analysis of Algorithm 1.
\begin{lemma} \label{lemma1}
	In Algorithm 1, the updated angle $ \theta $ always satisfies 
	\begin{align}\label{propo1}
		\underline{\theta} \leq \theta \leq \bar{\theta}.
	\end{align}
\end{lemma}
\begin{pf}
	According to Algorithm 1, either $ \theta = \arctan \frac{\phi\left(\boldsymbol{x}_{-}\right)-\phi\left(\boldsymbol{x}_{+}\right)}{g\left(\boldsymbol{x}_{+}\right)-g\left(\boldsymbol{x}_{-}\right)} $ or $ \theta= \frac{\bar{\theta}+\underline{\theta}}{2}$, where $ 0 \leq \underline{\theta} \leq \bar{\theta} \leq \frac{\pi}{2} $.
	Therefore, we only need to prove the first case, that is
	\begin{equation*}
		\tan\underline{\theta} \leq \frac{\phi\left(\boldsymbol{x}_{-}\right)-\phi\left(\boldsymbol{x}_{+}\right)}{g\left(\boldsymbol{x}_{+}\right)-g\left(\boldsymbol{x}_{-}\right)} \leq \tan\bar{\theta}.
	\end{equation*}
	As the point $ \boldsymbol{x}_{+} $ is obtained from the problem (\ref{sub1}) with $ \lambda=\tan\underline{\theta} $, it follows that
	\begin{align*}
		\phi(\boldsymbol{x}_{+})+&\tan\underline{\theta} \cdot  g(\boldsymbol{x}_{+}) \leq \phi(\boldsymbol{x}_{-})+\tan\underline{\theta}\cdot g(\boldsymbol{x}_{-})\\
		&\Longrightarrow \quad \tan \underline{\theta} \leq \frac{\phi\left(\boldsymbol{x}_{-}\right)-\phi\left(\boldsymbol{x}_{+}\right)}{g\left(\boldsymbol{x}_{+}\right)-g\left(\boldsymbol{x}_{-}\right)}.
	\end{align*}
	Furthermore, from the relationship  between $ \boldsymbol{x}_{-} $ and $ \bar{\theta} $, the similar conclusion  is obtained
	\begin{equation*}
		\quad \tan \bar{\theta} \geq \frac{\phi\left(\boldsymbol{x}_{-}\right)-\phi\left(\boldsymbol{x}_{+}\right)}{g\left(\boldsymbol{x}_{+}\right)-g\left(\boldsymbol{x}_{-}\right)}. \tag*{\qed}
	\end{equation*}  
\end{pf}

In Algorithm 1, the multiplier update strategy indicates that the interval $[\underline{\theta}, \bar{\theta}]$ shrinks by at least a factor of $(1 - \gamma)$. {According to Lemma \ref{lemma1}, furthermore, two properties of this algorithm are established:}
\begin{itemize}
	\item[(I)] At any iteration, the angle $ \theta^* $ corresponding to the  optimal Lagrange multiplier  always satisfies $ \underline{\theta} \leq \theta^* \leq \bar{\theta}  $.
	\item[(II)] After the k-th iteration, it follows that
	\begin{equation*}
		\bar{\theta}-\underline{\theta} \leq \frac{\pi}{2}(1-\gamma)^k.
	\end{equation*}
\end{itemize}

Based on these properties, we can proceed to derive the following theorem:
\begin{theorem}
	Given $ \epsilon >0 $, Algorithm 1 terminates within a finite number of iterations $ k_0 $ and $ k_0 \leq  \lceil{\log _{(1-\gamma)} \frac{2}{\pi} \varepsilon}\rceil$.
\end{theorem}
\begin{pf}
	This theorem can be derived from Lemma \ref{lemma1} and the property (I)(II).  \qed
\end{pf}

{We prove the optimality of the solution $\boldsymbol{x}^*$ obtained by Algorithm 1. This solution is derived from either line \ref{optimum1} or line \ref{optimum2}, corresponding to the two termination cases of this algorithm. For each case, we present the following optimality theorems.}

{\begin{theorem} \label{thm_op1}
		The solution $ \boldsymbol{x}^* $ obtained from the line \ref{optimum1} in Algorithm 1 is the optimal solution to problem (\ref{pro1}).
	\end{theorem}
\begin{pf}
	The KKT system of problem (\ref{pro1}) is
	\begin{subequations}
		\begin{align}
			&\phi_{i}^{\prime}+\lambda g_{i}^{\prime}-v_{i}+w_{i}=0, \quad i=1, \ldots, n, \label{KKTa}\\
			&{\sum_{i=1}^{n} g_{i}\left(x_{i}\right) \leqslant b},\label{KKTb} \\
			&{\lambda\left(\sum_{i=1}^{n} g_{i}\left(x_{i}\right)-b\right)=0}, \label{KKTc}\\
			&v_{i}\left(l_{i}-x_{i}\right)=0, \quad i=1, \ldots, n, \label{KKTd}\\
			&w_{i}\left(x_{i}-u_{i}\right)=0, \quad i=1, \ldots, n, \label{KKTe}\\
			&l_{i} \leqslant x_{i} \leqslant u_{i}, \quad i=1, \ldots, n, \label{KKTf}\\
			&{v_{i},w_{i} \geqslant 0, \quad i=1, \ldots, n}, \label{KKTg}\\
			&{\lambda \geqslant 0},\label{KKTh}
		\end{align}
	\end{subequations}
	where the vector $ \boldsymbol{v}=(v_1,\ldots,v_n), \boldsymbol{w}=(w_1,\ldots,w_n) $ denote the Lagrange multiplier associated with box constraints $ \boldsymbol{x}\in \boldsymbol{X} $.
	The solution $ \boldsymbol{x}^* $ obtained from the line \ref{optimum1} in Algorithm 1 satisfies $g(\boldsymbol{x}^*) = b$, and the Lagrange multiplier $ \lambda > 0 $ at this time. Therefore, current iterate $ (\boldsymbol{x}^*,\lambda)  $ satisfies (\ref{KKTb})(\ref{KKTc})(\ref{KKTh}).
	
	Additionally, the solution $ \boldsymbol{x}^*=\underset{\boldsymbol{x}\in \boldsymbol{X}}{\operatorname{argmin}\  }L(\boldsymbol{x},\lambda) $ and $ \lambda $ is obtained at the line \ref{lam_up1} or \ref{lam_up2}. Therefore, current iterate $ (\boldsymbol{x}^*,\lambda)  $ satisfies the remaining conditions in KKT system.
	
	Given that problem (\ref{pro1}) is convex and the solution $ (\boldsymbol{x}^*,\lambda)  $ satisfies all the conditions in KKT system, then $ \boldsymbol{x}^*$ is the optimal solution to problem (\ref{pro1}).  \qed
\end{pf}

\begin{lemma}\label{lemma2}
	For some positive $ \lambda^* $, if there are two optimal solutions $ \boldsymbol{x}_{+} $ and $ \boldsymbol{x}_{-} $ of the problem (\ref{sub1}) where $ g(\boldsymbol{x}_{+}) \geq b $ and $ g(\boldsymbol{x}_{-}) \leq b $, then $ \lambda^* $ is optimal for the Lagrangian dual problem (\ref{lag1}).
\end{lemma}
\begin{pf}
{	Let $ L(\lambda) = \underset{\boldsymbol{x}\in \boldsymbol{X}}{\operatorname{min}\  }L(\boldsymbol{x},\lambda) $, then the Lagrangian dual problem (\ref{lag1}) can be written as $ \underset{\lambda \geq 0}{\operatorname{max}} \  L(\lambda) $. Therefore, we just need to prove that  $ L(\lambda) \leq L(\lambda^*), \forall \lambda \neq \lambda^*  $. And the optimal solution to the problem (\ref{sub1}) is denoted by $ \boldsymbol{y} $ with respect to the multiplier $ \lambda $.} \\
	When $ \lambda < \lambda^* $, 
	\begin{align*}
		L(\lambda)=\phi(\boldsymbol{y})+\lambda(g(\boldsymbol{y})-b)&\leq\phi(\boldsymbol{x}_{+})+\lambda(g(\boldsymbol{x}_{+})-b)\\
		&\leq \phi(\boldsymbol{x}_{+})+\lambda^*(g(\boldsymbol{x}_{+})-b)\\
		&=L(\lambda^*).
	\end{align*}
	Similarly, in the case of $ \lambda > \lambda^* $,
	\begin{align*}
		L(\lambda)=\phi(\boldsymbol{y})+\lambda(g(\boldsymbol{y})-b)&\leq\phi(\boldsymbol{x}_{-})+\lambda(g(\boldsymbol{x}_{-})-b)\\
		&\leq \phi(\boldsymbol{x}_{-})+\lambda^*(g(\boldsymbol{x}_{-})-b)\\
		&=L(\lambda^*). \tag*{\qed}
	\end{align*}   
\end{pf}

	In other words, after the "While" loop (line \ref{while_begin}-\ref{while_end})  in  Algorithm 1, if we obtain the aforementioned two solutions corresponding to an identical  multiplier, then this multiplier is optimal for the Lagrangian dual problem (\ref{lag1}).
	
{\begin{lemma}\label{lem_omega}
	For problem (\ref{pro1}), the optimal solution set of $ \underset{\boldsymbol{x}\in \boldsymbol{X}}{\operatorname{min}\  }L(\boldsymbol{x},\lambda^*) $ is convex, where $ \lambda^* $ is the optimal Lagrange multiplier.
\end{lemma}
\begin{pf}
	Let 
	\begin{align*}
		\Omega^*:=&\{\boldsymbol{x}|\boldsymbol{x}=\underset{\boldsymbol{x}\in \boldsymbol{X}}{\operatorname{argmin}\  }L(\boldsymbol{x},\lambda^*)\} \\
		=&\{\boldsymbol{x}|\boldsymbol{x}=\underset{\boldsymbol{x}}{\operatorname{argmin}\  }L(\boldsymbol{x},\lambda^*)\ \cap \ \boldsymbol{x}\in \boldsymbol{X}\}
	\end{align*}
	be the optimal solution set of $ \underset{\boldsymbol{x}\in \boldsymbol{X}}{\operatorname{min}\  }L(\boldsymbol{x},\lambda^*) $.
	Since the functions $ \phi(\boldsymbol{x}),g(\boldsymbol{x}) $ are convex, the Lagrangian dual function $ L(\boldsymbol{x},\lambda^*) = \phi(\boldsymbol{x})+\lambda^*(g(\boldsymbol{x})-b)$ with $\ \lambda^*>0$ is also convex.
	Therefore, the optimal solution set $ \Omega^* $ is convex.  \qed
\end{pf}}
	Additionally, according to the convexity and  continuously differentiability of $ \phi(\boldsymbol{x}) $ and  $ g(\boldsymbol{x}) $, we can derive that within $ \boldsymbol{x} \in \Omega^* $, both  the functions $ \phi(\boldsymbol{x}) $ and  $ g(\boldsymbol{x}) $ are linear.  

\begin{theorem}\label{thm_op2}
	The solution $ \boldsymbol{x}^* $ obtained from the line \ref{optimum2} in Algorithm 1 is the optimal solution to problem (\ref{pro1}).
\end{theorem}
\begin{pf}
	If Algorithm 1 does not terminate at line \ref{optimum1}, two solutions $ \boldsymbol{x}_{+} $ and $ \boldsymbol{x}_{-} $ are obtained after the "While" loop (line \ref{while_begin}-\ref{while_end}), where $ g(\boldsymbol{x}_{+}) > b $ and $ g(\boldsymbol{x}_{-}) < b $. According to Lemma \ref{lemma2}, the current Lagrange multiplier $ \lambda^* $ is optimal for Lagrangian dual problem (\ref{lag1}) and $ \lambda^*>0 $.
	The dual optimal solution $ \boldsymbol{x}_{-} $ satisfies all the conditions in KKT system except (\ref{KKTc}).
	
	According to Lemma \ref{lem_omega}, the optimal solution set $ \Omega^* $ of problem (\ref{pro1}) is convex. Since $ g(\boldsymbol{x}) $ is convex and continuously differentiable, for all $ \boldsymbol{x} \in \Omega^* $, the function $ g_i(x_i) $ is linear for $ i = 1, \ldots, n $. Let\[ {\boldsymbol{x}^*}=\alpha \boldsymbol{x}_{+}+(1-\alpha)\boldsymbol{x}_{-}, \alpha \in \left(0,1\right), \] then $ {\boldsymbol{x}}^* \in \Omega^* $.
	We just need to solve the one-dimensional linear subproblem $ g({\boldsymbol{x}}^*)=b $, whose variable is $ \alpha $. Therefore, at the line \ref{optimum2} in Algorithm 1, the solution $ (\boldsymbol{x}^*,\lambda^*)  $  satisfies the condition (\ref{KKTc}) in KKT system.
	
	Given that problem (\ref{pro1}) is convex and the solution $ (\boldsymbol{x}^*,\lambda^*)  $ satisfies all the conditions in KKT system, then $ \boldsymbol{x}^*$ is the optimal solution to problem (\ref{pro1}).  \qed
\end{pf}
	
	For Algorithm 2, we can obtain the similar finite termination and convergence conclusions, whose proofs are omitted here.
\begin{corollary}
	Given $ \epsilon >0 $, Algorithm 2 terminates within a finite number of iterations $ k_0 $ and $ k_0 \leq  \lceil{\log _{(1-\gamma)} \frac{2}{\pi} \varepsilon}\rceil$.
\end{corollary}
{\begin{lemma}\label{lem_omega2}
		For problem (\ref{pro2}), the optimal solution set of $ \underset{\boldsymbol{x}\in \boldsymbol{X}}{\operatorname{min}\  }L(\boldsymbol{x},\lambda^*) $, \[ \Omega^*=\{\boldsymbol{x}|\boldsymbol{x}=\underset{\boldsymbol{x}\in \boldsymbol{X}}{\operatorname{argmin}\  }L(\boldsymbol{x},\lambda^*)\} \] is convex, where $ \lambda^* $ is the optimal Lagrange multiplier..
	\end{lemma}
	The proof is similar to Lemma \ref{lem_omega}, but the sign of optimal Lagrange multiplier $ \lambda^* $ is unrestricted. Because the function $ g(\boldsymbol{x}) $ is linear, the optimal solution set $ \Omega^* $ is still convex.

	\begin{corollary}
		The solution $ \boldsymbol{x}^* $ obtained from Algorithm 2 is the optimal solution to problem (\ref{pro2}).
	\end{corollary}
	
}

\section{Numerical experiments}  \label{sec5}
{In this section, we provide extensive numerical results for six types of practical test problems derived from various applications.  Specifically, the first two test problems, described in Sections \ref{com_ware} and \ref{quapro}, are associated with problem (\ref{pro1}), while the remaining four test problems, discussed in Sections \ref{port} through \ref{nega}, are associated with problem (\ref{pro2}). The experiments were conducted using MATLAB 9.0 on a computer equipped with an Intel Core i7-7800X processor (3.5 GHz, 12 cores) and 32 GB of RAM. Multiple sets of random numerical experiments were conducted for the same type of test problems, wherein we compared the CPU time  of all methods.} The CPU time is given  in seconds, obtained  by  calculating the geometric mean of all CPU times. And the CPU time is set to 7200 seconds if a method cannot solve the test problems within 2 hours. 

{In addition to  the pegging algorithm \citet{bretthauer2002nonlinear}, SLBQP1 algorithm (\citet{dai2006new}) and Aug\_lag algorithm (\citet{torrealba2022augmented}), we also utilized the Gurobi  solver (\citet{gurobi}) and CVX (\citet{gb08}, \citet{cvx}) to solve test problems. A brief overview of these  solvers is provided below:
\begin{itemize}
	\item Gurobi 11.0, a widely used commercial optimization solver that excels in solving linear, integer, and nonlinear optimization problems. Gurobi offers several efficient algorithms, including the Barrier Method and Dual Simplex Method, to ensure both the efficiency and accuracy of the solution process.
	\item CVX, a MATLAB-based package designed for solving convex optimization problems. CVX offers a user-friendly interface to model convex programs, supporting various problem types, including linear, quadratic, and semidefinite programming. 
	
\end{itemize}}

\subsection{Commodity warehousing problem} \label{com_ware}
The commodity warehousing test problem is as follows:
\begin{align*}
	\min _{x} & \sum_{i=1}^{n} c_{i} x_{i}+k_{i} / x_{i} \\
	\text { s.t. } & \sum_{i=1}^{n} a_{i} x_{i} \leq b, \\
	& x_{i}=\left[l_{i}, u_{i}\right], \quad i=1, \ldots, n.
\end{align*}	
In the numerical experiment, we set  $ a_{i} \in[1,4] $, $ c_{i} \in[10,30] $, $ k_{i} \in   [5,30], l_{i} \in(0,3]$   and  $ u_{i} \in(3,6]. $
\begin{table*}[b]
	\caption{The CPU time for commodity warehousing test problem.}
	\label{tab_q5}
	\begin{tabular}{ccccc}
		\toprule
		n        & Algorithm 1          & Pegging     & Gurobi      & CVX    \\
		\midrule
		$ 2\times 10^2 $ & \textbf{0.0048} & 0.8610  & 0.0842 & 0.8783    \\
		$ 2\times 10^3 $ & \textbf{0.0094} & 21.1969 & 0.4806 & 5.8177  \\
		$ 2\times 10^4 $ & \textbf{0.3362} & 7200    & 7200   & 814.4285 \\
		$ 2\times 10^5 $ & \textbf{0.5187} & 7200    & 7200   & -         \\
		$ 2\times 10^6 $ & \textbf{5.3691} & 7200    & 7200   & -   \\
		\bottomrule
	\end{tabular}
\end{table*}

Table \ref{tab_q5} presents the CPU time comparison of Algorithm 1, the pegging algorithm, Gurobi, and CVX. The results demonstrate that, for commodity warehousing test problems, the proposed Algorithm 1 exhibits the highest computational efficiency, outperforming Gurobi solver by at least an order of magnitude and surpassing both CVX and the pegging algorithm by over two orders of magnitude.
Notably, for problem size $ n=2 \times 10^4 $, both Gurobi solver and the pegging algorithm fail to solve the test problems within 2 hours. For $ n=2 \times 10^5 $, CVX encounters memory overflow and the corresponding CPU time is marked  with "-", while Algorithm 1 successfully identifies the optimal solution in under one second. Even for $ n=2\times 10^6 $, Algorithm 1 consistently obtains the optimal solution within a few seconds, highlighting its significant advantage in computational efficiency.

\begin{figure}[htbp]  
	\centering
	\includegraphics[scale=0.36]{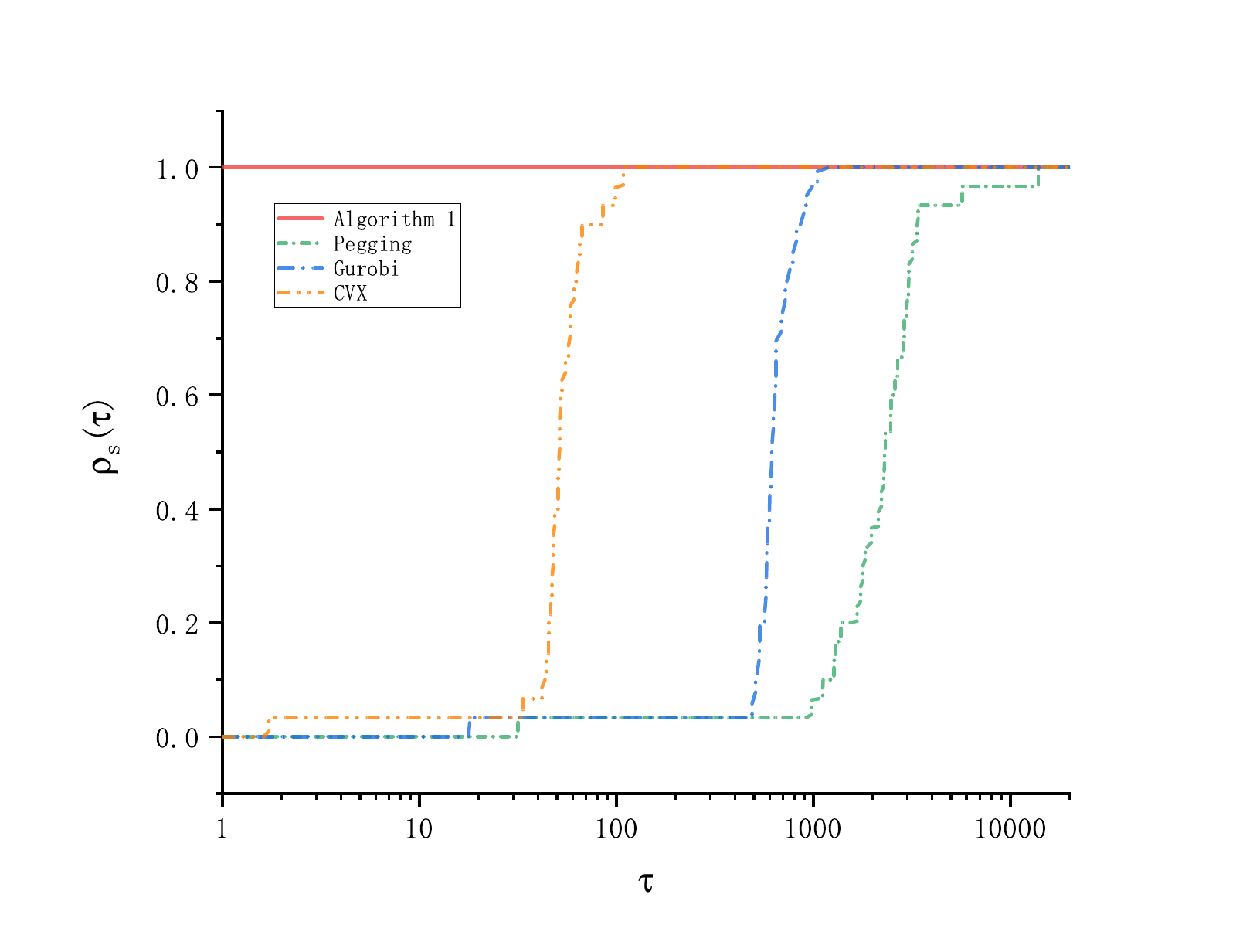}
	\caption{Perfomance profile of CPU time for $ 30 $ random commodity warehousing test problems with a problem size of $ n=2\times 10^3 $.}
	\label{fig_pro5}
\end{figure}

To compare the CPU time for  the above four methods more intuitively, we introduce the performance ratio (\citet{dolan2002benchmarking}) \begin{equation*}
	r_{p, s}^t=\frac{t_{p, s}}{\min \left\{t_{p, s}: s \in \mathcal{S}\right\}},
\end{equation*}
where $ t_{p, s} $  is the CPU time required to solve problem $ p $ by method $ s $ and $ S $ is the set of all the methods. And then 
\begin{equation*}
	\rho_{s}^t(\tau)=\frac{1}{n_{p}} \operatorname{size}\left\{p \in \mathcal{P}: r_{p, s}^t \leq \tau\right\},
\end{equation*}
where $ n_{p} $ is the  total number of test problems and the set $ \mathcal{P} $ includes all test problems. Furthermore, $ \rho_{s}^t(\tau) $  signifies the probability for method $  s \in \mathcal{S}  $ that a performance ratio $  r_{p, s}^t  $ is  the best ratio under a factor  $ \tau \in \mathbb{R} $. 

For the scenario where $ n=2\times 10^3 $, we randomly generate 30 test problems and evaluate their performance using the performance profile, as illustrated in Figure \ref{fig_pro5}. The results indicate that Algorithm 1 consistently demonstrates the highest computational efficiency across all test problems, followed by Gurobi solver, CVX, and pegging algorithm.

{\subsection{Quadratic problem with a quadratic constraint} \label{quapro}
In the context of problem (\ref{pro1}), we compare the performance of Algorithm 1 with Gurobi solver and CVX, as no relevant numerical experiments have been conducted yet.

	We present the following quadratic problem as the test problem, which is similar to the portfolio investment problem.
	\begin{align*}
		\underset{\boldsymbol{x}}{\operatorname{min}} \quad&\sum_{i=1}^{n} \frac{1}{2}d_{i} x_{i}^2 - c_{i} x_{i} \\
		s.t.\quad &\sum_{i=1}^{n} \frac{1}{2}a_{i} x_{i}^2 - z_{i} x_{i} \leq b,  \\
		&x_{i} \in \left[l_{i}, u_{i}\right], \quad i=1, \ldots, n. 
\end{align*}}
In the numerical experiments, we set $ a_{i} \in[1,30], z_{i} \in[1,35] ,d_{i} \in[1,20],$  $ c_{i} \in[1,25] ,  l_{i}\in[0,3] $  and  $ u_{i} \in(3,11] $.

\begin{table*}[htbp]
	\caption{The CPU time  for the quadratic test problem with a quadratic constraint.}
	\label{tab1}
	\begin{tabular}{cccc}
		\toprule
		n  & Algorithm 1 
		&  Gurobi   & CVX  \\
		\midrule
		$ 2\times 10^3 $ & \textbf{0.0064} & 1.2310    & \multicolumn{1}{l}{53.1443} \\
		$ 2\times 10^4 $ & \textbf{0.0936} & 46.9997   & -                           \\
		$ 2\times 10^5 $ & \textbf{0.3388} & 1826.0822 & -                           \\
		$ 2\times 10^6 $ & \textbf{4.5959} & 7200      & -                        \\  
		\bottomrule
	\end{tabular}
\end{table*}
We randomly generate 30 test problems for each problem size $ n $ and conduct numerical experiments to compare the performance of Algorithm 1, Gurobi solver and CVX. The CPU time comparison across varying problem sizes is summarized in Table \ref{tab1}.
The results clearly demonstrate that Algorithm 1 significantly outperforms both Gurobi solver and CVX in terms of computational efficiency, achieving a speedup of at least three orders of magnitude. Notably, for large-scale problems ($n= 2\times 10^6 $), Gurobi solver fails to obtain the optimal solution within 2 hours, whereas Algorithm 1 identifies the optimal solution within a few seconds. Furthermore, CVX is unable to solve test problems when $n = 2\times 10^4 $ due to memory limitations.
\begin{figure}[h]  
	\centering
	\includegraphics[scale=0.36]{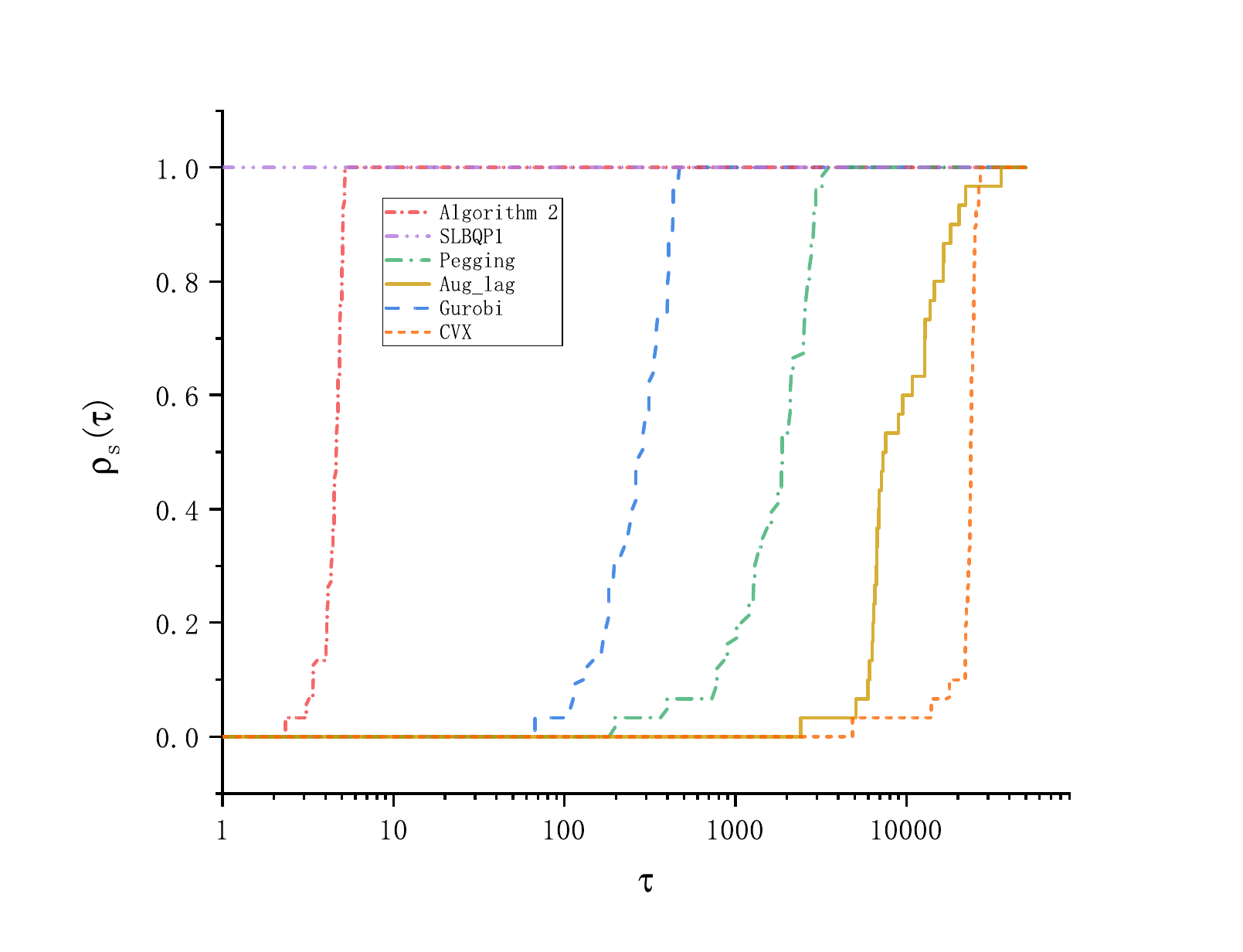}
	\caption{Perfomance profile of CPU time for $ 30 $ random portfolio investment test problems with a problem size of $ n=2\times 10^4 $.} 
	\label{fig_pro1}
\end{figure}

\subsection{Portfolio investment problem}	\label{port}
In order to facilitate a more comprehensive comparison of  different methods, we extend the portfolio investment problem into a more general model as follows:
\begin{table*}[hp]  
	\caption{The CPU time for portfolio investment test problem.}
	\label{tab_q1}
	\begin{tabular}{ccccccc}
		\toprule
		n & Algorithm 2  & SLBQP1     & Pegging  & Aug\_lag &  Gurobi & CVX \\
		\midrule
		$ 2\times 10^3 $ & 0.0025 & \textbf{0.0017} & 1.3710   & 3.1153  & 0.2276    & 6.0995   \\
		$ 2\times 10^4 $ & 0.0209 & \textbf{0.0051} & 7.1116   & 76.1399 & 1.0285    & 106.7068 \\
		$ 2\times 10^5 $ & 0.2880 & \textbf{0.1039} & 241.7251 & 7200    & 22.0351   & -        \\
		$ 2\times 10^6 $ & 3.4228 & \textbf{1.0316} & 7200     & 7200    & 2597.1532 & -\\
		\bottomrule
	\end{tabular}
\end{table*}
\begin{align*} 
	\underset{\boldsymbol{x}}{\operatorname{min}}\quad&\sum_{i=1}^{n} \frac{1}{2}d_{i} x_{i}^2 - c_{i} x_{i} \\
	s.t.\quad&\sum_{i=1}^{n} a_{i} x_{i} = b,  \\
	&x_{i} \in \left[l_{i}, u_{i}\right], \quad i=1, \ldots, n. 
\end{align*}
In the numerical experiment, we set  $ a_{i} \in[1,30], d_{i} \in[1,20], c_{i} \in[1,25] ,  l_{i}\in[0,3] $  and  $ u_{i} \in(3,11] $.
The CPU time comparison of all methods is summarized in Table \ref{tab_q1}. The results indicate that portfolio investment test problems, the SLBQP1 algorithm demonstrates the highest computational efficiency, followed by Algorithm 2, with a performance gap of at most one order of magnitude between the two. Moreover, both algorithms significantly outperform the remaining three methods, achieving a speedup of over two orders of magnitude.


For the scenario where $ n=2\times 10^4 $, we randomly generate 30 test problems, and the corresponding performance profile is shown in Figure \ref{fig_pro1}. The results indicate that SLBQP1 algorithm achieves the highest computational efficiency, followed by Algorithm 2,  Gurobi solver, the pegging algorithm, the Aug\_lag algorithm and CVX.

\subsection{Sampling problem} \label{samp}
The sampling test problem is as follows:
{\begin{align*}
		\underset{\boldsymbol{x}}{\operatorname{min}} \quad &\sum_{i=1}^{n} c_{i}/x_{i}  \\
		s.t.\quad &\sum_{i=1}^{n} a_{i}x_{i} = b , \\ 
		&x_{i} =\left[l_{i}, u_{i}\right], \quad i=1, \ldots, n .
\end{align*}}
In the numerical experiment, we set  $ a_{i} \in[1,4], c_{i} \in[5,30] ,  l_{i}\in[0,3] $  and  $ u_{i} \in(3,6] $.

 The comparison of CPU time  for sampling test problems is summarized in Table \ref{tab_s1}.
\begin{table*}[h]
	\caption{The CPU time for sampling test problems.}
	\label{tab_s1}
	\begin{tabular}{cccccc}
		\toprule
		n  & Algorithm 2     & Pegging  & Aug\_lag &  Gurobi  &  CVX \\
		\midrule
		$ 2\times 10^2 $ & \textbf{0.0012} & 2.3409                & {7.3260}   & 0.1505                   & 0.7637                \\
		$ 2\times 10^3 $ & \textbf{0.0042} & 15.3010               & {282.7052} & 3.2343                   & 5.2448                \\
		$ 2\times 10^4 $ & \textbf{0.0437} & 251.3642              & 7200                         & 39.4006                  & 97.5100               \\
		$ 2\times 10^5 $ & \textbf{0.4601} & 7200 & 7200                         & 507.0256                 & \multicolumn{1}{c}{-} \\
		$ 2\times 10^6 $ & \textbf{5.0926} & 7200 & 7200                         & \multicolumn{1}{c}{7200} & \multicolumn{1}{c}{-} \\
		\bottomrule
	\end{tabular}
\end{table*}
The results reveal that, for sampling test problems, Algorithm 2 achieves the highest computational efficiency, surpassing Gurobi solver and CVX by at least two orders of magnitude and outperforming the pegging algorithm and the Aug\_lag algorithm by over three orders of magnitude. 
For the case where $ n= 2\times 10^4 $, we randomly generate 30 test problems and depict the corresponding perfrmance profile in Figure \ref{fig_pro2}. This figure illustrates that Algorithm 2 consistently exhibits the best computational efficiency across all test problems, followed by Gurobi solver, pegging algorithm and CVX. Notably, the Aug\_lag algorithm fails to solve the test problems within 2 hours.

\begin{figure}[h]  
	\centering
	\includegraphics[scale=0.36]{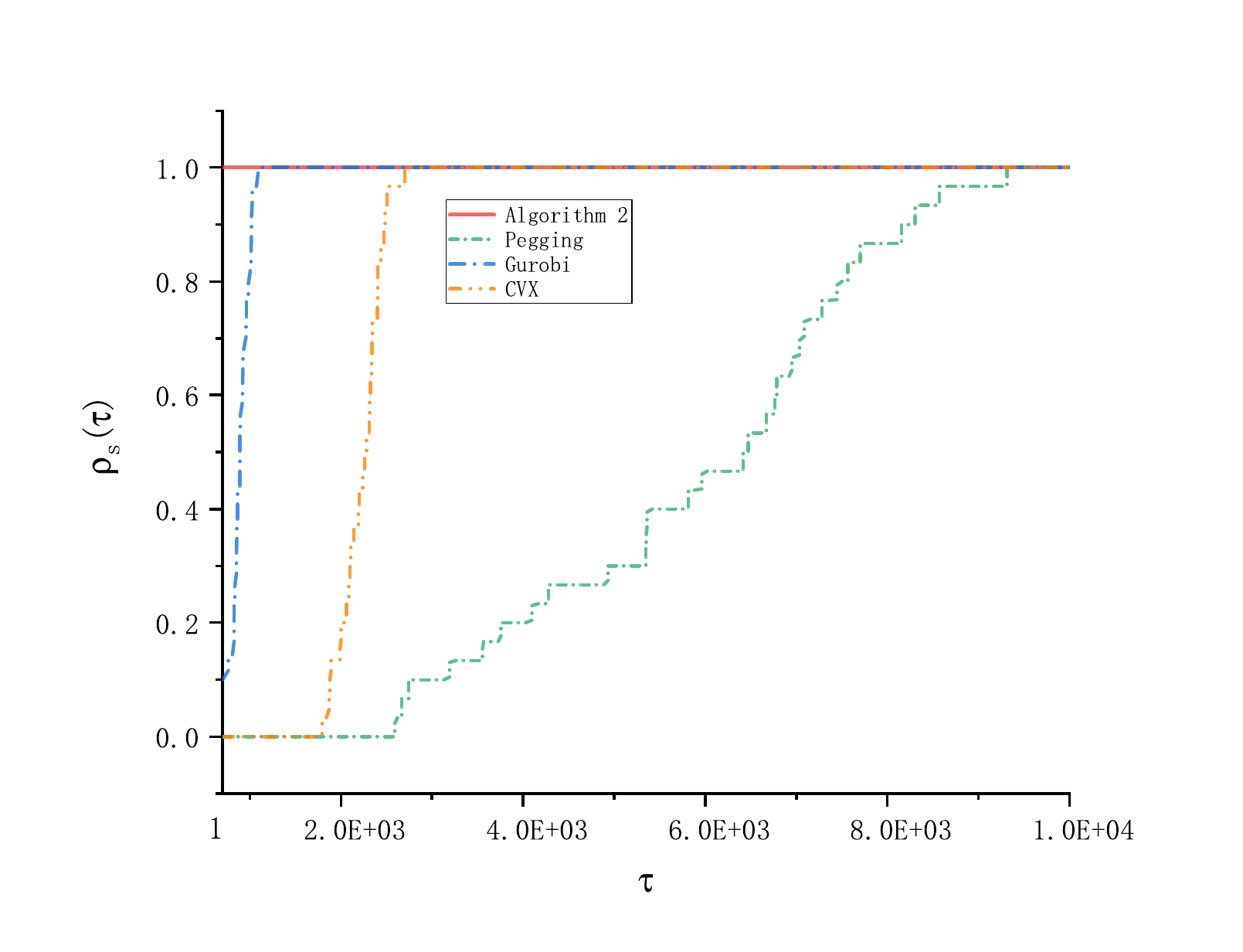}
	\caption{ Perfomance profile of CPU time for $ 30 $ random sampling test problems with a problem size of $ n=2\times 10^4 $.}  
	\label{fig_pro2}
\end{figure}

{\subsection{Target search problem} \label{theory}
	The target search test problem is as follows:
	\begin{align*}
		\underset{\boldsymbol{x}}{\operatorname{min}} \quad &\sum_{i=1}^{n} m_{i}\left(e^{-c_{i} x_{i}}-1\right)  \\
		s.t.\quad &\sum_{i=1}^{n} a_{i}x_{i} = b , \\ 
		&x_{i} =\left[l_{i}, u_{i}\right], \quad i=1, \ldots, n .
\end{align*}}
In the numerical experiment, we set  $ a_{i} \in[1,3], m_{i} \in[0.5,8], c_{i} \in[0.1,3] ,  l_{i}\in[0,0.1] $  and  $ u_{i} \in(0.1,5] $. 

\begin{table*}[h]
	\caption{The CPU time for  target search test problems.}
	\label{tab_t}
	\begin{tabular}{cccc}
		\toprule
		n        & Algorithm 2     & Aug\_lag  & Gurobi    \\
		\midrule
		$ 2\times 10^2 $ & \textbf{0.0043} & 257.6770 & 0.0838    \\
		$ 2\times 10^3 $ & \textbf{0.0060} & 7200     & 0.5254    \\
		$ 2\times 10^4 $ & \textbf{0.0489} & 7200     & 7200 \\
		$ 2\times 10^5 $ & \textbf{0.5569} & 7200     & 7200      \\
		$ 2\times 10^6 $ & \textbf{6.0177} & 7200     & 7200      \\
		\bottomrule
	\end{tabular}
\end{table*}

Table \ref{tab_t} compares the  CPU time  of  Algorithm 2, the Aug\_lag algorithm and Gurobi solver. The results indicate that Algorithm 2 significantly outperforms the other two algorithms. Notably, for $ n=2 \times 10^4 $, neither the Aug\_lag algorithm nor Gurobi solver  can solve the test problems within 2 hours, whereas Algorithm 2 efficiently identifies the optimal solution within one second. Even for  $ n=2\times 10^6 $, Algorithm 2 obtains the optimal solution within a few seconds, highlighting its significant advantage in computational efficiency.

{\subsection{Negative entropy problem} \label{nega}
	The negative entropy problem is mentioned in \citet{nielsen1992massively} and we consider the following model:
	\begin{align*}
		\underset{\boldsymbol{x}}{\operatorname{min}} \quad &\sum_{i=1}^{n} x_{i} \log \left(\frac{x_{i}}{a_{i}}-1\right)  \\
		s.t.\quad &\sum_{i=1}^{n} a_{i}x_{i} = b , \\ 
		&x_{i} =\left[l_{i}, u_{i}\right], \quad i=1, \ldots, n .
\end{align*}}
In the numerical experiment, we set  $ a_{i} \in[1,3], c_{i} \in[0.1,1.9], l_{i}\in[2,10] $  and  $ u_{i} \in(10,21] $.


\begin{table*}[hb]
	\caption{The CPU time  for negative entropy test problems.}
	\label{tab_n}
	\begin{tabular}{cccc}
		\toprule
		n        & Algorithm 2     & Aug\_lag  & Gurobi    \\
		\midrule
		$ 2\times 10^2 $ & \textbf{0.0029}  & 55.8015 & 7200   \\
		$ 2\times 10^3 $ & \textbf{0.0274}  & 7200    & 7200   \\
		$ 2\times 10^4 $ & \textbf{0.2411}  & 7200    & 7200   \\
		$ 2\times 10^5 $ & \textbf{3.1193}  & 7200    & 7200   \\
		$ 2\times 10^6 $ & \textbf{31.5301} & 7200    & 7200  \\
		\bottomrule
	\end{tabular}
\end{table*}

Table \ref{tab_n} compares the CPU time of Algorithm 2, the Aug\_lag algorithm and Gurobi solver. The results show that Algorithm 2 outperforms the other two algorithms by at least four orders of magnitude in computational efficiency.
Notably, for $n= 2\times 10^3 $, neither the Aug\_lag algorithm nor Gurobi solver can solve the test problems within 2 hours, whereas Algorithm 2 efficiently obtains the optimal solution within tens of milliseconds. Even for $ n= 2\times 10^6 $, Algorithm 2 identifies the optimal solution within tens of seconds, demonstrating its significant advantage in computational efficiency.

\section{Conclusion}  \label{sec6}
The CONRAP has undergone thorough studies; however, most existing algorithms are unable to solve this problem without relying on monotonicity assumptions for the objective and constraint functions. Based on the problem's convexity, we propose a novel Lagrange multiplier update strategy for CONRAP, which utilizes the values of both the objective and constraint functions at the current and previous iterations. This strategy accelerates the process of identifying optimal solutions.
Subsequently,  leveraging the problem's separability, the Lagrangian dual problem is decomposed into  $ n $ one-dimensional tractable subproblems; these subproblems are readily solved by using the function's mononicity and gradient.  
Furthermore, utilizing the Lagrange duality theory, we propose two Lagrangian dual algorithms  to solve two types of CONRAP. Through rigorous theoretical analysis, we demonstrate that both of the proposed algorithms converge to optimal solutions within a finite number of iterations.

To evaluate the performance of our proposed algorithms, we conducted extensive numerical experiments on six types of practical test problems from various applications. 
The results indicate that Algorithm 1 can solve two types of test problems, associated with problem (\ref{pro1}), with the highest computational efficiency, being at least an order of magnitude faster than the pegging algorithm, Gurobi solver and CVX. Furthermore, Algorithm 1 exhibits significant performance in handling large-scale test problems ($ n= 2\times 10^6 $), identifing the optimal solution within a few seconds.
When solving problem (\ref{pro2}) with a quadratic objective function, Algorithm 2 demonstrates a slightly slower computational efficiency than SLBQP1 algorithm, with an order of magnitude gap at most. Nevertheless, Algorithm 2 significantly outperforms the pegging algorithm, the Aug\_lag algorithm, Gurobi solver and CVX, being at least two orders of magnitude faster than these methods. 
Furthermore, when solving the remaining three types of test problems, Algorithm 2 exhibits the highest computational efficiency, outperforming other algorithms and solvers by at least two orders of magnitude. Notably, for large-scale test problems ($ n= 2\times 10^6 $), Algorithm 2 efficiently identifies the optimal solution within seconds or tens of seconds, whereas other  algorithms and solvers fail to solve these problems within the 2-hour time limit.


\appendix
\section{Projected gradient algorithm} \label{appa}
As described in Section \ref{alstep}, when $  g\left(\boldsymbol{x}_{g}\right)=b  $ holds, we need to identify the optimal solution to $ \phi(\boldsymbol{x}) $ under the set $ \Omega_g $.
If $ \|\nabla g\left(\boldsymbol{x}_{g}\right)\| = 0 $, $ \Omega_g $ may  have multiple elements and we adopt the following projected gradient algorithm proposed by \citet{iusem2003convergence} in order to identify the solution $ \boldsymbol{x} $ satisfying $ \boldsymbol{x}= \underset{\boldsymbol{x}\in \Omega_g}{\operatorname{argmin}\  }\phi(\boldsymbol{x})  $ with the convergence guarantee. The projected gradient algorithm is  presented as Algorithm \ref{pg}. In this algorithm, $  \beta_{k}, \gamma_{k} $  are positive stepsizes and $ P_{\Omega_g} $ denotes the orthogonal projection onto the set $ \Omega_g $. 
%
\begin{algorithm}[h]
	\caption{Projected gradient algorithm} \label{pg} 
	$ \text { Procedure }\left[ \boldsymbol{x}^*\right]=\operatorname{Projected \_Gradient}(\phi(\boldsymbol{x}),g(\boldsymbol{x}) , \boldsymbol{x}^{0},  \boldsymbol{X}) $ 
	\begin{algorithmic}[1]
		\State Set $k=0, \epsilon> 0,\beta_{k}> 0,\gamma_{k} > 0 $ and $ \tilde{\boldsymbol{x}} = 2\boldsymbol{u} $;
		\While{$ \|\tilde{\boldsymbol{x}} - \boldsymbol{x}^{k}\|>\epsilon $}
		\State $ \tilde{\boldsymbol{x}} = \boldsymbol{x}^{k} $;
		\State $ \boldsymbol{x}_{p} = P_{\Omega_g}\left(\boldsymbol{x}^{k} - \beta_{k} \nabla \phi(\boldsymbol{x}^{k}) \right)$;
		\State $ \boldsymbol{x}^{k+1} = \boldsymbol{x}^{k} - \gamma_{k} \left(\boldsymbol{x}_{p}-\boldsymbol{x}^{k}\right) $;
		\State $ k=k+1 $;
		\EndWhile
		\State Return $ \boldsymbol{x}^{k} $.
	\end{algorithmic}
\end{algorithm}

Armijo search along the feasible direction:  $ \left\{\beta_{k}\right\} \subset[\tilde{\beta}, \hat{\beta}] $  for some $  0<   \tilde{\beta} \leq \hat{\beta} $  and  $ \gamma_{k} $  determined with an Armijo rule, namely
\begin{equation*}
\gamma_{k} = 2^{-\ell(k)}
\end{equation*}
with
\begin{align*}
\ell(k) = \min &\big\{j \in \mathbb{Z}_{\geq 0}:  f\left(\boldsymbol{x}^{k} - 2^{-j}\left(\boldsymbol{x}_{p} - \boldsymbol{x}^{k}\right)\right)  \\
& \quad \leq f\left(\boldsymbol{x}^{k}\right) - \sigma 2^{-j} \nabla f\left(\boldsymbol{x}^{k}\right)^{\top} \left(\boldsymbol{x}^{k} - \boldsymbol{x}_{p}\right) \big\}
\end{align*}

for some  $ \sigma \in(0,1)  $.

%


\bibliographystyle{cas-model2-names}

\bibliography{reference}



\end{document}